\documentclass[12pt,psamsfonts,leqno,oneside,letterpaper]{amsart}
\usepackage[dvips,text={6.5truein,9truein},left=1truein,top=1truein]{geometry}
\usepackage{amssymb,amsmath,amscd,enumerate}
\usepackage{amsrefs}
\usepackage[pdftex]{graphicx}
\usepackage[colorlinks,linkcolor=blue,citecolor=blue,pdfstartview=FitH]{hyperref}
\input xy
\xyoption{all}
\SelectTips{cm}{12}

\theoremstyle{definition}
\newtheorem{para}{}[section]

\newtheorem{remark}[para]{Remark}
\newtheorem{remarks}[para]{Remarks}
\newtheorem{notation}[para]{Notation}
\newtheorem{convention}[para]{Convention}
\newtheorem{definition}[para]{Definition}
\newtheorem{definitions}[para]{Definitions}

\newcommand\Alternatives{\begin{enumerate}[(i)]}
\newcommand\EndAlternatives{\end{enumerate}}
\newcommand\Conditions{\begin{enumerate}[(1)]}
\newcommand\EndConditions{\end{enumerate}}

\theoremstyle{plain}
\newtheorem{theorem}[para]{Theorem}
\newtheorem{lemma}[para]{Lemma}
\newtheorem{proposition}[para]{Proposition}
\newtheorem{corollary}[para]{Corollary}
\newtheorem{conjecture}[para]{Conjecture}
\newtheorem{claim}[equation]{}

\numberwithin{equation}{para}
\numberwithin{figure}{section}

\newcommand\Number{\begin{para}}
\newcommand\EndNumber{\end{para}}
\newcommand\Definition{\begin{definition}}
\newcommand\EndDefinition{\end{definition}}
\newcommand\Definitions{\begin{definitions}}
\newcommand\EndDefinitions{\end{definitions}}
\newcommand\Theorem{\begin{theorem}}
\newcommand\EndTheorem{\end{theorem}}
\newcommand\Conjecture{\begin{conjecture}}
\newcommand\EndConjecture{\end{conjecture}}
\newcommand\Remark{\begin{remark}}
\newcommand\EndRemark{\end{remark}}
\newcommand\Remarks{\begin{remarks}}
\newcommand\EndRemarks{\end{remarks}}
\newcommand\Convention{\begin{convention}}
\newcommand\EndConvention{\end{convention}}
\newcommand\Notation{\begin{notation}}
\newcommand\EndNotation{\end{notation}}
\newcommand\Lemma{\begin{lemma}}
\newcommand\EndLemma{\end{lemma}}
\newcommand\Proposition{\begin{proposition}}
\newcommand\EndProposition{\end{proposition}}
\newcommand\Corollary{\begin{corollary}}
\newcommand\EndCorollary{\end{corollary}}
\newcommand\Claim{\begin{claim}}
\newcommand\EndClaim{\end{claim}}
\newcommand\Proof{\begin{proof}}
\newcommand\EndProof{\end{proof}}
\newcommand\Equation{\begin{equation}}
\newcommand\EndEquation{\end{equation}}

\newcommand\Bullets{\begin{itemize}}
\newcommand\EndBullets{\end{itemize}}

\usepackage{color}

\newcommand\ZZ{{\mathbb Z}}

\newcommand\QQ{{\mathbb Q}}
\newcommand\HH{{\mathbb H}}
\newcommand\EE{{\mathbb E}}
\newcommand\RR{{\mathbb R}}

\newcommand\calm{{\mathcal M}}
\newcommand\dist{\mathop{\rm dist}}
\newcommand\inter{\mathop{\rm int}}

\newcommand\vol{\mathop{\rm vol}}

\newcommand\isomplus{\mathop{{\rm Isom}_+}}
\newcommand\diam{\mathop{\rm diam}}
\newcommand\arccosh{\mathop{\rm arccosh}}

\newcommand\Fix{{\rm Fix}}
\newcommand\Per{{\rm Per}}

\newcommand\tM{\widetilde M}

\newcommand\tF{\widetilde F}
\newcommand\tPhi{\widetilde \Phi}

\newcommand\tc{\widetilde c}

\newcommand\Xsqrt[1]{\sqrt{\rule{0pt}{2ex}#1}}

\begin{document}
\author{Marc Culler}
\address{Department of Mathematics, Statistics, and Computer Science (M/C 249)\\
University of Illinois at Chicago\\
851 S. Morgan St.\\
Chicago, IL 60607-7045}
\email{culler@math.uic.edu}

\author{Peter B. Shalen}
\address{Department of Mathematics, Statistics, and Computer Science (M/C 249)\\  University of Illinois at Chicago\\
  851 S. Morgan St.\\
  Chicago, IL 60607-7045} 
\email{shalen@math.uic.edu}
\thanks{Both authors are partially supported by NSF grant DMS-0906155}

\title{Margulis numbers for Haken manifolds}

\begin{abstract}
  For every closed orientable hyperbolic Haken $3$-manifold
  and, more generally, for any orientable hyperbolic
  $3$-manifold $M$ which is homeomorphic to the interior of a Haken
  manifold, the number $0.286$ is a Margulis number. If
  $H_1(M;\QQ)\ne0$, or if $M$ is closed and contains a semi-fiber,
  then $0.292$ is a Margulis number for $M$.
\end{abstract}

\maketitle

\section{Introduction}

If $M$ is an orientable hyperbolic $n$-manifold, we may write
$M=\HH^n/\Gamma$ where $\Gamma\le\isomplus(\HH^n)$ is discrete and
torsion-free. The group $\Gamma$ is determined, up to conjugacy in
$\isomplus(\HH^n)$, by the hyperbolic structure on $M$.

For any $\gamma\in\Gamma$ and any $P\in\HH^3$ we shall write
$d_P(\gamma)=\dist(P,\gamma\cdot P)$. (Here and throughout the paper,
$\dist$ denotes the hyperbolic distance on $\HH^n$.)  We will define a
{\it Margulis number} for $M$, or for $\Gamma$, to be a positive real
number $\mu$ with the following property:

\Claim\label{general claim}
For any $P\in \HH^n$, the subgroup  of $\Gamma$ generated by $\{x\in\Gamma\;|\;d_P(x) < \mu\}$
has an abelian subgroup of finite index.
\EndClaim

If $n\le3$ then any subgroup of $\Gamma$ which contains an abelian subgroup of finite index is itself abelian. Hence in this case \ref{general claim} may be replaced by the following simpler condition:

\Claim
For any $P\in \HH^n$, if $x$ and $y$ are elements of
$\Gamma$ such that $\max(d_P(x),d_P(y))<\mu$, then $x$ and
$y$ commute.
\EndClaim

The Margulis Lemma \cite{bp}*{Chapter D} implies that for every
$n\ge2$ there is a positive constant which is a Margulis number for
every hyperbolic $n$-manifold. The largest such number,
$\mu(n)$, is called the {\it Margulis constant} for hyperbolic
$n$-manifolds.

Margulis numbers play a central role in the geometry of hyperbolic
manifolds. If $M$ is a hyperbolic $n$-manifold, which for simplicity
we take to be closed, and $\mu$ is a Margulis number for $M$, then the
points of $M$ where the injectivity radius is less than $\mu/2$ form a
disjoint union of ``tubes'' about closed geodesics whose geometric
structure can be precisely described. Topologically they are open
$(n-1)$-ball bundles over $S^1$. This observation and the Margulis
Lemma can be used to show, for example, that for every $V>0$ there is
a finite collection of compact orientable $3$-manifolds
$M_1,\ldots,M_N$, whose boundary components are tori, such that every
closed, orientable hyperbolic $3$-manifold of volume at most $V$ can
be obtained by a Dehn filling of one of the $M_i$.

The value of $\mu(3)$ is not known; the best known lower bound is
$0.104\ldots$, due to Meyerhoff \cite{meyerhoff}.  In this paper we
will derive a larger lower bound for the Margulis numbers of closed
orientable hyperbolic Haken $3$-manifolds.  In fact, our
results apply to any orientable $3$-manifold which
is homeomorphic to the interior of a Haken $3$-manifold.

Recall that a compact orientable irreducible $3$-manifold $N$ is {\it
  Haken} if it contains an incompressible surface or if it is
homeomorphic to $B^3$. We will say that $N$ is {\it strictly Haken} if
it contains an incompressible surface which is not a fiber or a
semi-fiber.  (We will review the definitions of incompressible
surface, fiber and semi-fiber in Section \ref{surface section}.)

Our main result is the following.

\Theorem\label{main theorem or wall theorem} Let $M$ be an
  orientable hyperbolic $3$-manifold which is homeomorphic to the
interior of a Haken manifold.  (In particular $M$ may be a closed
Haken manifold.)  Then $0.286$ is a Margulis number for $M$. If
$H_1(M;\QQ)\ne0$, or if $M$ is closed and contains a semi-fiber, then
$0.292$ is a Margulis number for $M$.
\EndTheorem

\Remark The condition $H_1(M;\QQ)\ne0$ always holds unless $M$ is
closed or $M=\HH^3$. This is because $B^3$ is the only non-closed
Haken manifold with first Betti number $0$.
\EndRemark

To prove Theorem \ref{main theorem or wall theorem} for a hyperbolic
manifold $M=\HH^3/\Gamma$ which is homeomorphic to the interior of a
given Haken manifold (or strict Haken manifold) $N$, we must
obtain a lower bound for $\max(d_P(x),d_P(y))$ whenever $P$ is a point
of $\HH^3$ and $x$ and $y$ are non-commuting elements of $\Gamma$.

If we assume for simplicity that $N$ contains a
separating incompressible surface, then the choice of such a surface
allows us to identify $\Gamma\cong\pi_1(M)$ with an amalgamated free
product $A\star_CB$. The proof of the lower bound for
$\max(d_P(x),d_P(y))$ breaks up into various cases depending on the
normal forms for $x$ and $y$ in $A\star_CB$. In some cases one can
show that suitable short words in $x$ and $y$ generate a free group of
rank $2$, and one can then obtain the desired bound by combining the
``$\log3$ Theorem'' and its refinements (\cite{paradoxical},
\cite{accs}, \cite{agol}, \cite{cg}, \cite{acs}) with some delicate
hyperbolic trigonometry; in other cases one
can show that certain short words in $x$ and $y$ generate a free
{\it semigroup} of rank $2$, and one can then obtain 
a bound using packing arguments of the type used in \cite{sw}*{Section 5}.

In cases where the $\log3$ Theorem and packing arguments do not give a
suitable lower bound, a third method is required. We shall illustrate
this method by describing it in a simple case. Suppose that we have
$x\in A$, $y\in B$, and that neither $x^2$ nor $y^2$ lies in $C$. Let
$X$ (respectively $Y$) denote the set of all elements of
$\Gamma=A\star_CB$ whose normal form begins with an element of $A$
(respectively $B$). Then $\Gamma$ is the disjoint union of $X$, $Y$
and $C$, and we have (i) $x^{\pm 1} Y \subset X$ and $y^{\pm 1} X
\subset Y$, (ii) $x Y \cap x^{-1} Y = \emptyset$, and (iii) $y X \cap
y^{-1} X = \emptyset$. The argument now proceeds in a way that is
partially analogous to the proof of the $\log3$ theorem given in
\cite{paradoxical} for the special case of a rank-$2$ free group for
which the normalized area measure is a Patterson-Sullivan measure on
the sphere at infinity $S_\infty$. The decomposition $\Gamma=X\cup
Y\cup C$ gives rise to a decomposition of the area measure on
$S_\infty$, and the set-theoretic conditions (i)---(iii) give
information about how the terms in the decomposition transform under
$x$ and $y$. This information is used to obtain a lower bound for
$\max(d_P(x),d_P(y))$; here the hyperbolic trigonometry is
  combined with an argument analogous to the one used in
\cite{paradoxical}.

The transition from information about measures on $S_\infty$ to
information about the quantity $\max(d_P(x),d_P(y))$ uses Lemma 5.5 of
\cite{paradoxical}, a result that was also applied in \cite{accs}. The
statement of the lemma given in \cite{paradoxical} contains some minor
errors and an irrelevant hypothesis. We provide a corrected and
improved statement, with a complete proof, in Section \ref{measure
  section} of the present paper, along with an account of the small
additional arguments needed for the applications of the lemma given in
\cite{paradoxical} and \cite{accs}.

In the sketch given above of a special case of the argument involving
decompositions of measures, the group-theoretical argument was
stated in terms of free products with amalgamation.  In the general
case, it is best to use the language of group actions on trees. The
needed background concerning incompressible surfaces and actions on
trees is given in Sections \ref{tree section} and \ref{surface
  section}. The key argument involving decompositions of measures is
given in section \ref{decomposition section}. The packing arguments
that apply in the case where short words in $x$ and $y$ generate a
free semigroup are given in Section \ref{packing section}.
The hyperbolic trigonometry needed for the arguments described
above appears in Section \ref{sphere section}.
In Section \ref{putting it together}, these ingredients are assembled
to prove the main theorem.

We are grateful to David Futer for pointing out that our methods apply
to the case of a non-compact hyperbolic manifold, and to Mark Kidwell
and David Krebes for pointing out some errors in the original
manuscript.  We thank the anonymous referee for several helpful
suggestions which improved the exposition.
  
\section{Measures and displacements}\label{measure section}

The purpose of this section is to prove the following result, which is
a corrected and improved version of Lemma 5.5 of \cite{paradoxical}
and will be needed in Section \ref{decomposition section} of the
present paper.

\Lemma\label{elles sans foot} Let $a$ and $b$ be numbers in $[0,1]$
which are not both equal to $0$ and are not both equal to $1$.  Let
$\gamma $ be a loxodromic isometry of $\HH^{3}$ and let $z$ be a point
in $\HH^{3}$.  Suppose that $\nu $ is a measure on $S_{\infty }$ such
that
\begin{enumerate}
\item$\nu  \le  A_{z}$,
\item$\nu (S_{\infty }) \le  a$, and
\item $\int^{}_{S_{\infty }}\lambda^2_{\gamma,z}\, d\nu\ge b$. 
\end{enumerate}
Then we have $a>0$, $b<1$, and
$$\dist(z,\gamma \cdot z) \ge \frac12\log\frac{b(1-a)}{a(1-b)}.$$   
\EndLemma

\bigskip

In \cite{paradoxical}*{Lemma 5.5}, the numerator and denominator were
interchanged in the conclusion. The issue of guaranteeing that the
denominator in the corrected expression is non-zero was not addressed
in the statement or proof, and the hypothesis which assures that this
denominator is non-zero was omitted. Furthermore, the statement of
\cite{paradoxical}*{Lemma 5.5} contained the hypothesis $a\le1/2$,
which is not needed for the proof.

For the applications of \cite{paradoxical}*{Lemma 5.5} given in
\cite{paradoxical} and \cite{accs}, we will explain after the proof of
Lemma \ref{elles sans foot} how the latter result can be applied to
replace \cite{paradoxical}*{Lemma 5.5}.

\Proof[Proof of Lemma \ref{elles sans foot}] Throughout this proof we
use the conventions of \cite{paradoxical}.  We let $h$ denote the
constant dist$(z,\gamma \cdot z)$, and set $c = \cosh h$ and $s =
\sinh h$.  Since $\gamma$ is loxodromic we have $h>0$, so that
$c>s$. We let $\lambda $ denote the function $\lambda _{\gamma
  ,z}$. We identify $\overline\HH^3$ conformally with the unit ball in
$\RR^{3}$ in such a way that $z$ is the origin (so that $S_{\infty
}$ has the round metric centered at $z)$ and $\gamma ^{-1}\cdot z$ is
on the positive vertical axis.

According to \cite{paradoxical}*{2.4}, we have $\lambda (\zeta ) =
{\mathcal P}(z,\gamma ^{-1}\cdot z,\zeta )$ for all $\zeta \in
S_{\infty }$.  Hence by \cite{paradoxical}*{2.1.1}, $\lambda $ is
given by the formula
$$\lambda (\zeta ) = (c - s \cos  \phi )^{-1}$$
where $\phi = \phi (\zeta )$ is the angle between the positive
vertical axis and the ray from the origin through $\zeta $; thus $\phi
$ is the polar angle of $\zeta $ in spherical coordinates.

Set $A = A_{z}$.  Since $S_{\infty }$ has the round metric
centered at $z$, the measure  $A$
is obtained by dividing the area measure on the unit sphere by the 
area $4\pi $
of the sphere.  In spherical coordinates $\theta $ and $\phi $ on
the unit sphere we  have
$dA = (1/4\pi )\sin\phi\, d\phi d\theta $.

Set $\phi _{0} = \arccos  (1 - 2a)\in[0,\pi]$,  and  let $C \subset 
S_{\infty }$  denote  the  spherical  cap
defined by the inequality $\phi  \ge  \phi _{0}$.  Then we have
$$A(C) = \frac{1}{4\pi}\int^{2\pi}_0\int^{\phi_0}_0\sin\phi\, d\phi d\theta
= \frac12 (1 - \cos \phi _{0}) = a .$$ Thus by hypothesis (ii) we
have $A(C) \ge \nu (S_{\infty })$.  Observe also that since $\lambda $
is given by the function $(c - s \cos \phi )^{-1}$, which is positive
and monotone decreasing for $0 \le \phi \le \pi $, we have $\inf
\lambda (C) \ge \sup \lambda (S_\infty - C)$.  Since we also have $\nu
\le A$ by hypothesis (i), we may apply \cite{paradoxical}*{Lemma 5.4}
with $f = \lambda ^{2}$ to obtain
$$\begin{aligned}
\int^{}_{S_{\infty }}\lambda ^{2} d\nu  &\le  \int^{}_{C} \lambda
^{2} dA = \frac{1}{4\pi}\int^{2\pi }_{0} \int^{{\phi}_0 }_{0}
\frac{\sin  \phi}{(c - s \cos  \phi )^{2}} d\phi  d\theta \\
&= \frac12\int^{{\phi}_0 }_{0} \frac{\sin  \phi}{(c - s \cos 
\phi )^{2}} d\phi  = 
\frac{1}{2s} \left( \frac{1}{c - s} - \frac{1}{c - s\cos\phi _{0}}\right),
\end{aligned}
$$
where  the  last  step  follows  from  the  substitution $u = c -
s \cos  \phi $.
Recalling that $\cos  \phi _{0} = 1 - 2a$ and using hypothesis
(iii), we find that
$$b \le  \int^{}_{S_{\infty }}\lambda ^{2} d\nu  \le 
\frac{a}{(c - s)(c - s + 2as)},$$    
which by the definitions of $c$ and $s$ gives
\Equation\label{let bigons be bigons}
be^{-2h}+ab-abe^{-2h}\le a.
\EndEquation
It follows from (\ref{let bigons be bigons}) that if $a=0$ then
$b\le0$, and that if $b=1$ then $a\ge1$; in view of the hypotheses it
follows that $a>0$ and that $b<1$. It now follows from (\ref{let
  bigons be bigons}) that
$$e^{-2h} \le  \frac{a (1 - b)}{b (1 - a)},$$
which implies the conclusion of the lemma.
\EndProof

In the applications of \cite{paradoxical}*{Lemma 5.5} given in
\cite{paradoxical}, Lemma \ref{elles sans foot} can be applied without
change. In particular, these applications involved only specific
values of $a$ and $b$ which are not equal to $1$ or $0$.

For the application of \cite{paradoxical}*{Lemma 5.5} given in
\cite{accs}, citing Lemma \ref{elles sans foot} requires a bit more
care. The application appears in the proof of Theorem 6.1 of
\cite{accs}, and in the notation of that proof we have $a=\alpha_i$
and $b=1-\beta_i$ for some given $i\in\{1,\ldots,k\}$. In order to
apply Lemma \ref{elles sans foot} we must check that we cannot have
$a=b=1$ or $a=b=0$. If $a=b=1$ then $\alpha_i=1$ and $\beta_i=0$,
which is impossible because $\alpha_i\le\beta_i$ in the context of the
proof. Now suppose that $a=b=0$, so that $\alpha_i=0$ and
$\beta_i=1$. Since $\sum_{i=1}^k(\alpha_i+\beta_i)=1$ and $k\ge2$, we
may fix an index $j\ne i$ in $\{1,\ldots,k\}$, and we have
$\alpha_j=\beta_j=0$. Hence the measures $\nu_j$ and $\nu_j'$ are
identically zero. But we have
$$\int  (\lambda _{\xi_j^{-1} ,z_{0}}) ^{2} d\nu'_j  = 1 -
\int d\nu _j ,$$
which is impossible if $\nu_j=\nu_j'=0$. 

\section{Group actions on trees}\label{tree section}

We begin by reviewing some elementary notions related to
group actions on trees.  Our point of view here is similar
to that taken in \cite{nonsep}.

\Number\label{tree stuff}
By a {\it tree} we will mean a $1$-connected $1$-dimensional
simplicial complex. We may regard the set of vertices of a tree as an
integer metric space by defining the distance between two vertices
to be the number of edges in the arc joining them.

A {\it line} in a tree $T$ is a subcomplex isomorphic to the real line
with the standard triangulation, in which the vertices are the
integers.

If $\Gamma$ is a group, we will define a {\it $\Gamma$-tree} to be a
tree equipped with a simplicial action of $\Gamma$ which has no
inversions in the sense of \cite{Serre}. We will be using basic facts
about $\Gamma$-trees proved in \cite{Serre}. In particular, if $T$ is
a $\Gamma$-tree and $\gamma$ is a non-trivial element of $\Gamma$ then
either $\gamma$ has a non-empty fixed tree, in which case it is said
to be {\it $T$-elliptic}, or $\gamma$ has a unique invariant line,
called its {\it axis}; in this case $\gamma$ is said to be {\it
  $T$-hyperbolic}.  The stabilizer of an edge $e$ or a vertex $s$
will be denoted $\Gamma_e$ or $\Gamma_s$ respectively.

If $T$ is a $\Gamma$-tree and $\gamma$ is a $T$-elliptic element of
$\Gamma$, we will denote its fixed tree by $\Fix(\gamma)$.  We also
set $\Per(\gamma)=\bigcup_{n=1}^\infty\Fix(\gamma^n)$. Rewriting
$\Per(\gamma)$ as $\bigcup_{n=1}^\infty\Fix(\gamma^{n!})$ shows that
it is a monotone union of subtrees. Hence $\Per(\gamma)$ is a subtree
of $T$, which we refer to as the {\it periodic subtree} of $\gamma$.
\EndNumber

\Definition Let $\Gamma$ be a group and let $T$ be a $\Gamma$-tree. We
shall say that the action of $\Gamma$ on $T$ is {\it linewise
  faithful} if for every line $L$ in $T$, the subgroup of $\Gamma$
that fixes $L$ pointwise is trivial. We shall say that the action is
{\it trivial} if some vertex of $T$ is fixed by the entire group
$\Gamma$.
\EndDefinition

\Proposition\label{dum da-dum dum} Let $\Gamma$ be a group and let $T$
be a $\Gamma$-tree. Suppose that the action of $\Gamma$ on $T$ is
non-trivial. Then for each vertex $s$ of $T$, the stabilizer $\Gamma_s$
has infinite index in $\Gamma$.
\EndProposition

\Proof If $\Gamma_s$ has finite index in $\Gamma$ then the
orbit $\Gamma\cdot s$ is finite. Let $T_0$ denote the smallest subtree
containing $\Gamma\cdot s$. Then $T_0$ is a finite subtree and is
$\Gamma$ invariant. It is clear that the action of $\Gamma$ on $T_0$
factors through a finite quotient $G$ of $\Gamma$. But it follows from
\cite{Serre}*{No. I.6.3.1} that if $G$ is any finite group, then for
any $G$-tree $T_0$ the action of $G$ on $T_0$ is trivial. This
contradicts the non-triviality of the action of $\Gamma$ on $T$.
\EndProof

\Definition We will say that elements $x_1, \ldots, x_n$ of a group
$\Gamma$ are {\it independent} if they generate a free
group of rank $n$, and {\it semi-independent} if they generate
a free semigroup of rank $n$.  (In other words, $x_1, \ldots, x_n$
are semi-independent if distinct positive words in these elements
represent distinct elements of $\Gamma$.)
\EndDefinition

\Lemma\label{dey shoulda dunnit} Let $\Gamma$ be a group and let $T$
be a $\Gamma$-tree. Suppose that $x$ and $y$ are $T$-hyperbolic
elements of $\Gamma$ whose axes are distinct. Then $x$ and $y$ are
semi-independent.
\EndLemma

\Proof
This follows from \cite{BH}*{p. 687, Proof of Lemma}.
\EndProof

\Proposition\label{ILGWU} Let $\Gamma$ be a group, and let $T$ be a
$\Gamma$-tree such that the action of $\Gamma$ on $T$ is non-trivial
and linewise faithful. Then any two non-commuting $T$-hyperbolic
elements of $\Gamma $ are semi-independent.
\EndProposition

\Proof Suppose that $x_1,x_2\in\Gamma$ are non-commuting
$T$-hyperbolic elements. Let $A_i$ denote the axis of $x_i$ in
$T$. Then $x_i$ acts by a translation on $A_i$. If $A_1=A_2$, it
follows that $x_1x_2x_1^{-1}x_2^{-1}$ fixes the line $A_1=A_2$
pointwise. Since the action of $\Gamma$ on $T$ is linewise faithful,
it follows that $x_1x_2x_1^{-1}x_2^{-1}=1$, a contradiction. Hence
$A_1\ne A_2$. It now follows from Lemma \ref{dey shoulda dunnit} that
$x$ and $y$ are semi-independent.
\EndProof

\section{Incompressible surfaces and actions on trees}\label{surface section}

Here, and in the sequel, we will often suppress base points when
denoting fundamental groups of connected spaces.  (We note that, while
inclusion homomorphisms are only defined up to conjugacy if we have
not chosen a base point, the injectivity of an inclusion homomorphism
is independent of such a choice.)  In addition, we will often assume
that an identification of the fundamental group of a connected
$3$-manifold $M$ with the group of deck transformations of its
universal cover $\tM$ has been fixed.

\Definition Let $M$ be a connected orientable $3$-manifold and $F$ an
orientable $2$-manifold embedded in $M$.  We will say that $F$ is {\it
  bicollared} if there exists an embedding $c:F\times[-1,1]\to M$,
proper in the sense of general topology, such that $c(x, 0) = x$ for
all $x\in F$, and $c((\inter F)\times[-1,1])\subset\inter M$.  The map
$c$ will be called a {\it bicollaring} of $F$.
\EndDefinition

\Number\label{hawaii} Let $F$ be a bicollared surface in a connected,
orientable $3$-manifold $M$. Let $c:F\times[-1,1]\to M$ be a
bicollaring of $F$. If $p:\tM\to M$ denotes the universal covering of
$M$, then $\tF=p^{-1}(F)$ inherits a bicollaring $\tc$. We may
partition $\tM$ into the components of $\tM-\tc(\tF\times(-1,1))$ and
the sets of the form $\tc(\tPhi\times\{t\})$, where $\tPhi$ ranges
over the components of $\tF$ and $t$ ranges over the interval
$(-1,1)$. Let $T$ denote the decomposition space defined by this
partition, and let $q:\tM\to T$ denote the decomposition map. Then $T$
may be regarded as a tree in such a way that the vertices of $T$ are
images under $q$ of the components of $\tM-\tc(\tF\times(-1,1))$, and
the edges of $T$ are images under $q$ of the components of $\tF$. If
$\Gamma$ denotes the group of deck transformations of $\tM$, the
action of $\Gamma$ on $\tM$ induces a simplicial action, without
inversions, on $T$. Thus $T$ has the structure of a $\Gamma$-tree in a
natural way; we call it the {\it dual tree} of the bicollared surface
$F$.  \EndNumber

\Definition\label{incompressible} Let $M$ be a compact orientable
$3$-manifold, possibly with non-empty boundary.  By an {\it
  incompressible surface} in $M$ we mean a connected bicollared
surface $F\subset M$ such that
\begin{itemize}
\item the inclusion homomorphism $\pi_1(F)\to\pi_1(M)$ is injective;
\item $F$ is not the boundary of a $3$-ball in $M$; and
\item $F$ is not parallel to a subsurface of $\partial M$.
\end{itemize}
\EndDefinition

\Proposition\label{nontrivial tree} Suppose that $F$ is an
incompressible surface in a compact orientable
$3$-manifold $M$.  Let $\Gamma$ denote the group of deck
transformations of $\tM$, and $T$ the dual $\Gamma$-tree of $F$.  Then
$T$ is a non-trivial $\Gamma$-tree.
\EndProposition

\Proof
We identify $\Gamma$ with $\pi_1(M)$.  The stabilizer of a
vertex of $T$ is then a conjugate of the image under inclusion of
$\pi_1(V)$ for some component $V$ of $M-F$.

If the surface $F$ is nonseparating then there is a surjective
homomorphism $\pi_1(M)\to\ZZ$ whose kernel contains the
image under inclusion of $\pi_1(M-F)$.  This shows that
the vertex stabilizers are proper subgroups of $\Gamma$ in this case.

Otherwise, we may write $M = A\cup B$ where $A$ and $B$ are compact
connected submanifolds of $M$ with $A\cap B = F$.  Since $F$ is
incompressible, the inclusion homomorphisms $\pi_1(F)\to\pi_1(A)$ and
$\pi_1(F)\to\pi_1(B)$ are injective. By Van Kampen's theorem,
$\pi_1(M)$ may be identified with the free product of $\pi_1(A)$ and
$\pi_1(B)$ amalgamated along $\pi_1(F)$. In particular $\pi_1(A)$ and
$\pi_1(B)$ are then identified with subgroups of $\pi_1(M)$. We will
show that both of these factors are proper subgroups.

Assume for a contradiction that $\pi_1(A) = \pi_1(M)$.  The normal
form theorem for an amalgamated free product then implies that
$\pi_1(F) = \pi_1(B)$, i.e. that the inclusion $F\hookrightarrow B$
induces an isomorphism of fundamental groups.  It follows from
\cite{hempel}*{Theorem 10.5} that there is a homeomorphism from $B$ to
$F\times[0,1]$ sending $F$ to $F\times\{0\}$.  But this implies that
$F$ is parallel to a subsurface of $\partial M$, a contradiction to
the incompressibility of $F$.  Thus $\pi_1(A)$ is a proper subgroup of
$\pi_1(M)$.  A symmetrical argument shows that $\pi_1(B)$ is a proper
subgroup of $\pi_1(M)$.

Since any vertex stabilizer is conjugate to one of the two proper
subgroups $\pi_1(A)$ and $\pi_1(B)$, it follows that $T$ is a
non-trivial $\Gamma$-tree in this case.
\EndProof

\Definition\label{fibers} An incompressible surface $F$ in a closed
orientable $3$-manifold $M$ is called a {\it fiber} if there is a
fibration $M \to S^1$ having $F$ as one of the fibers.  We will
say that $F$ is a {\it semi-fiber} if $F$ separates $M$ into two
components, each of which is the interior of a twisted $I$-bundle over
a non-orientable surface.
\EndDefinition

\Proposition\label{yes georgia} Let $M = \HH^3/\Gamma$ be a closed
hyperbolic $3$-manifold containing an incompressible surface $F$ which
is not a fiber or a semi-fiber. Let $g$ denote the genus of $F$, and
let $T$ denote the dual $\Gamma$-tree of $F$. Then for every
non-trivial $T$-elliptic element $\gamma\in\Gamma$, the diameter (as
an integer metric space, cf. \ref{tree stuff}) of the set of fixed
vertices $\gamma$ in $T$ is at most $14g-12$. In particular, the
action of $\Gamma$ on $T$ is linewise faithful.
\EndProposition

\Proof The first assertion follows from \cite{DeB}*{Corollary
 1.5}. The second assertion follows from the first, since a
line in a tree has infinite diameter.
\EndProof

\Lemma\label{georgia} Let $N$ be a connected orientable $3$-manifold
without boundary, and let $A_0$ and $A_1$ be disjoint incompressible
open annuli in $N$. Suppose that the inclusion homomorphism
$\pi_1(A_0)\to\pi_1(N)$ is an isomorphism. Then the inclusion
homomorphism $\pi_1(A_1)\to\pi_1(N)$ is also an isomorphism.
\EndLemma

\Proof In this proof, unlabeled homomorphisms will be understood to be
induced by inclusion.  For $i=0,1$ it follows from incompressibility
that the image of $\pi_1(A_i)\to\pi_1(N)$ has finite index in
$\pi_1(N)$. If one of the (bicollared) annuli $A_i$ did not separate
$N$ then the image of $H_1(A_i;\ZZ)\to H_1(N;\ZZ)$ would lie in the
kernel of a homomorphism of $H_1(N;\ZZ)$ onto $\ZZ$, a
contradiction. Hence each of the $A_i$ separates $N$, and there is a
connected submanifold $N_0$ of $N$ whose frontier is $A_0\cup
A_1$. Since the $A_i$ are incompressible in $N$,
$\pi_1(N_0)\to\pi_1(N)$ is injective. Since $\pi_1(A_0)\to\pi_1(N)$ is
surjective, $\pi_1(N_0)\to\pi_1(N)$ is also surjective, and is
therefore an isomorphism. It follows that $\pi_1(A_0)\to\pi_1(N_0)$ is
an isomorphism, and that the image of $\pi_1(A_1)\to\pi_1(N_0)$ has
finite index in the infinite cyclic group $\pi_1(N_0)$. Hence there is
a map $f:S^1\times[0,1]\to N_0$ such that for $i=0,1$ we have
$f(S^1\times\{i\})\subset A_i$, and $f|S^1\times\{i\}$ is
homotopically non-trivial. It now follows from Waldhausen's
generalized loop theorem \cite{genloop} that there is a properly
embedded planar surface $P\subset N_0$ having at most one boundary
component in each $A_i$, and such that each boundary component of $P$
is homotopically non-trivial in $\partial N_0$. Since the $A_i$ are
incompressible in $N$, the surface $P$ has no disk components; hence
it has exactly one boundary component in each $A_i$, and is therefore
an annulus. The boundary of $P$ must consist of a core curve in $A_0$
and a core curve in $A_1$. Since $\pi_1(A_0)\to\pi_1(N)$ is an
isomorphism, it now follows that $\pi_1(A_1)\to\pi_1(N)$ is also an
isomorphism.
\EndProof

\Remark
A proof of Lemma \ref{georgia} could also be based on
\cite{McCullough}*{Theorem 2}.
\EndRemark

\Proposition \label{alabama} Let $F$ be an incompressible surface in
an orientable $3$-manifold $M$.  Let $\Gamma$ denote the group of deck
transformations of the universal cover of $M$, and let $T$ denote the dual
$\Gamma$-tree of $F$. Suppose that $\gamma$ is an infinite-order
element of $\Gamma$ such that $\Fix(\gamma)\subset T$ contains at
least one edge of $T$. Then for every integer $n>0$ we have
$\Fix(\gamma^n)=\Fix(\gamma)$.
\EndProposition

\Proof We choose an open edge $e_0$ in $\Fix(\gamma)$. If $n>0$ is
given, $\Fix(\gamma^n)$ is a subtree containing $e_0$, and is
therefore a union of closed edges. Hence it suffices to show that if
an open edge $e_1$ is fixed by $\gamma^n$ then it is fixed by
$\gamma$. We may assume that $e_1\ne e_0$.

Let $C$ denote the infinite cyclic group $\langle\gamma\rangle$, and
for $i=0,1$ let $C_i$ denote the stabilizer of $e_i$ in $C$; then
$C_0=C$ since $\gamma$ fixes $e_0$, and $C_1$ is also infinite cyclic
since $\gamma^n$ fixes $e_1$. We need to show that $C_1=C$. Since $C$
stabilizes $e_0$ and $e_1\ne e_0$, the $C$-orbits of $e_0$ and $e_1$
are distinct.

We use the notation of Subsection \ref{hawaii}. For $i=0,1$, we have
$e_i=q(\tPhi_i\times(-1,1))$ for some component $\tPhi_i$ of $\tF$. The
stabilizer of $\tPhi_i$ is $C_i$. Hence $A_i:=\tPhi_i/C_i$ is
identified with a surface in the $3$-manifold $N:=\tM/C$,
which is incompressible since $F$ is incompressible in
$M$. Since $C_i\cong\pi_1(A_i)$ is infinite cyclic, each
$A_i$ is an annulus. Since $C_0=C$, the inclusion
homomorphism $\pi_1(A_0)\to\pi_1(N)$ is surjective. Since the $e_i$
are in distinct $C$-orbits we have $A_0\cap A_1=\emptyset$. It
therefore follows from Lemma \ref{georgia} that the inclusion
homomorphism $\pi_1(A_1)\to\pi_1(N)$ is surjective. This implies that
$C_1=C$, as required.
\EndProof

\Proposition \label{alaska} Let $\Gamma$ be a group, let $T$ be a
$\Gamma$-tree, and let $\gamma_0,\gamma_1$ be $T$-elliptic elements of
$\Gamma$. Suppose that
$\Per(\gamma_0)\cap\Per(\gamma_1)=\emptyset$.  Then
$\gamma_0$ and $\gamma_1$ are independent in $\Gamma$.
\EndProposition

\Proof The hypotheses immediately imply that $\gamma_0$ and $\gamma_1$
have infinite order. We will apply the Klein criterion, (or
``ping-pong lemma'') as stated in \cite{lyndonschupp}*{Proposition
  12.2}, taking the groups $G_1$ and $G_2$ of the latter result to be
the infinite cyclic groups generated by $\gamma_0$ and $\gamma_1$
respectively.  According to this criterion, we need only construct
disjoint subsets $\Omega_0$ and $\Omega_1$ of $T$ such that
$\gamma_i^n\cdot\Omega_i \subset \Omega_{1-i}$ for all $0\ne n\in\ZZ$.
For $i=0,1$ set $X_i=\Per(\gamma_i)$. By hypothesis we have $X_0\cap
X_1=\emptyset$.  Let $C$ be an open topological arc in $T$ whose
endpoints are vertices, one of which lies in $X_0$ and one in
$X_1$. Among all such arcs we choose $C$ in such a way as to minimize
the number of open edges that it contains. Then for $i=0,1$, the set
$C\cap X_i$ consists of a single vertex $s_i$. Let $e_i$ denote the
open edge which is contained in $C$ and has $s_i$ as an
endpoint. (Note that $e_0$ and $e_1$ may or may not coincide.)  We
define $\Omega_{1-i}$ to be the component of $T-e_i$ that contains
$s_i$. Since $T$ is a tree, $\Omega_i$ is a component of $T-C$. In
particular, the sets $\Omega_0$, $\Omega_1$ and $C$ are pairwise
disjoint. Since $X_i$ contains $s_i$ and is disjoint from $e_i$, we
have $X_i\subset \Omega_{i-1}$ for $i=0,1$. Note that $\Omega_i\cup C$
is connected for $i=0,1$.

Since $e_i$ is not contained in $\Per(\gamma_i)$ it follows that
$\gamma_i^n\cdot e_i\ne e_i$ for $0\ne n\in\ZZ$; thus we have
$\gamma_i^n\cdot\overline e_i\subset T- e_i$. On the other hand we
have $s_i\in X_i = \Per(\gamma_i)$ and hence $\gamma^n\cdot s_i\in
X_i\cap \gamma^n\cdot\overline e_i$.  In particular, $X_i\cup
\gamma^n\cdot\overline e_i$ is connected, and it follows that the
component of $ T- e_i$ that contains $\gamma_i^n\cdot\overline e_i$
also contains $s_i$ and is therefore equal to
$\Omega_{i-1}$. Thus we have shown that \Equation\label{arkansas}
\gamma_i^n\cdot\overline e_i\subset \Omega_{1-i}
\EndEquation
whenever $0\ne n\in\ZZ$ and $i\in\{0,1\}$.

It follows that $e_i$ is disjoint from $\gamma_i^{n}\cdot(C\cup
\Omega_i)$, and hence from its closure $\gamma_i^{n}\cdot(\overline
C\cup \Omega_i)$. Since $\overline C\cup \Omega_i$ is connected,
$\gamma_i^{n}\cdot(\overline C\cup \Omega_i)$ must be contained in a
component of $ T- e_i$; and since $s_i=\gamma_i\cdot s_i\in
\gamma_i^{n}\cdot\overline C$, this component must be
$\Omega_{1-i}$. Thus we have shown that $\gamma_i^n\cdot (\overline
C\cup \Omega_{1-i})\subset \Omega_i$ whenever $0\ne n\in\ZZ$ and
$i\in\{0,1\}$.  In particular, we have $\gamma_i^n\cdot
\Omega_i\subset \Omega_{1-i}$, as required.
\EndProof

\Proposition\label{old mexico} Let $M=\HH^3/\Gamma$ be a closed
orientable hyperbolic $3$-manifold, and let $T$ be a
non-trivial $\Gamma$-tree. Suppose that $\gamma_0$ and $\gamma_1$ are
non-commuting $T$-elliptic elements of $\Gamma$ such that
$\Fix(\gamma_0)\cap \Fix(\gamma_1)\ne\emptyset$. Then $\gamma_0$ and
$\gamma_1$ are independent in $\Gamma$.
\EndProposition

\Proof Choose a vertex $s\in \Fix(\gamma_0)\cap\Fix(\gamma_1)$. Then
$H:=\langle\gamma_1,\gamma_2\rangle$ is a subgroup of $\Gamma_s$.
Since $T$ is a non-trivial $\Gamma$-tree, it follows from
\ref{nontrivial tree} that $H < \Gamma_s$ has infinite
index in $\Gamma$.  But since $M$ is a closed hyperbolic $3$-manifold,
it follows from \cite{JS}*{Theorem VI.4.1} that every two-generator
subgroup of infinite index in $\Gamma=\pi_1(M)$ is free of rank at
most $2$; since $\gamma_0$ and $\gamma_1$ do not commute, $H$ is a
free group of rank $2$, i.e. $\gamma_0$ and $\gamma_1$ are
independent.
\EndProof

\Proposition \label{arizona} Let $F$ be an incompressible surface in the
closed orientable hyperbolic $3$-manifold $M=\HH^3/\Gamma$, and let
$T$ be the dual $\Gamma$-tree of $F$. Suppose that $\gamma_0$ and
$\gamma_1$ are non-commuting elements of $\Gamma$ such that
$\Fix(\gamma_i)\subset T$ contains at least one edge of $T$ for
$i=0,1$. Then $\gamma_0$ and $\gamma_1$ are independent in $\Gamma$.
\EndProposition

\Proof Since $M$ is a hyperbolic manifold, $\Gamma$ is torsion-free.
Thus it follows from Proposition \ref{alabama} that
$\Per(\gamma_i)=\Fix(\gamma_i)$ for $i=0,1$. Hence if
$\Fix(\gamma_0)\cap\Fix(\gamma_1)=\emptyset$, it follows from
Proposition \ref{alaska} that $\gamma_0$ and $\gamma_1$ are
independent. If $\Fix(\gamma_0)\cap\Fix(\gamma_1)\ne\emptyset$, then
$\gamma_0$ and $\gamma_1$ are independent by Proposition \ref{old
  mexico}.
\EndProof

\section{Word growth and displacements}\label{packing section}

We observe that for any $P\in\HH^3$ and any isometries $x$ and $y$ of
$\HH^3$, we have $d_P(xy)\le d_P(x)+d_P(y)$, $d_P(x^{-1})=d_P(x)$, and
$\dist(x\cdot P,y\cdot P)=d_P(x^{-1}y)$. These facts will be used
frequently in this and the following sections.

\Definition If $S$ is a finite subset of a group $\Gamma$, and $m$ is
a positive integer, we will denote by $b_m(S)$ the number of elements
of $\Gamma$ that can be expressed as words of length at most $m$ in
elements of $S$. We will set
$\omega(S)=\lim_{m\to\infty}b_m(S)^{1/m}$; it is pointed out in
\cite{sw} that this limit exists.
\EndDefinition

\Proposition\label{s'n'w prunes} Suppose that $\Gamma$ is a
torsion-free discrete subgroup of $\isomplus(\HH^n)$ for some integer
$n\ge2$. Let $x$ and $y$ be elements of $\Gamma$. Then for any point
$P\in\HH^n$ we have
$$\max(d_P(x),d_P(y))\ge\frac{\log\omega(\{x,y\})}{n-1}.$$
\EndProposition

\Proof (Cf. \cite{sw}*{proof of Proposition 5.2}) Set
$\omega=\omega(\{x,y\})$.  If the conclusion is false, there exist
real numbers $C<\omega$ and $\lambda < (\log C)/(n- 1)$ such that
$d_P(x) < \lambda$ and $d_P(y) < \lambda $. Since
$C<\omega$, the definition of $\omega(\{x,y\})$ implies that $b_m(S) >
C^m$ for all sufficiently large $m$.  Choose a bounded neighborhood
$U$ of $P$ such that $\gamma \cdot U \cap U = \emptyset$ for every
$\gamma \in \Gamma - { 1 } $. Set $v = \vol U>0$. If $\gamma \in
\Gamma $ is expressible as a word of length $\le m$ in $x$ and $y$
then $d_P(\gamma) < m\lambda $. Thus for any $m$ there are
$b_m(\{x,y\}) $ elements $\gamma \in \Gamma $ such that $\gamma \cdot
P$ lies in the ball of radius $m\lambda$ about $P $. Hence if $\Delta
= \diam U$, there are $b_m(S)$ disjoint sets of the form $\gamma \cdot
U$, $\gamma \in \Gamma $, in a ball of radius $m\lambda + \Delta$
about $P$.  But there is a constant $K$ depending on the dimension $n$
such that the volume of any ball of sufficiently large radius $r$ in
hyperbolic $n$-space is bounded above by $K\exp ((n-1)r )$. Hence for
large $m$ we have $ b_m(S) \cdot v < K\exp((n - 1)(m\lambda +
\Delta))$, and therefore \Equation\label{delaware} v\cdot C^{m} < K
\cdot \exp((n - 1)(m\lambda + \Delta)).
\EndEquation 
Now taking logarithms of both sides of (\ref{delaware}), dividing by
$m$ and taking limits as $m\to\infty$ we get $(n - 1)\lambda \ge\log
C$. This contradicts our choices of $C$ and $\lambda$.
\EndProof

\Corollary\label{florida} Suppose that $\Gamma$ is a torsion-free
discrete subgroup of $\isomplus(\HH^n)$ for some integer $n\ge2$. Let
$x$ and $y$ semi-independent elements of $\Gamma$. Then for any point
$P\in\HH^n$ we have
$$\max(d_P(x),d_P(y))\ge\frac{\log2}{n-1}.$$
\EndCorollary

\Proof For any positive integer $m$ there are $2^m$ positive words in
$x$ and $y$, and since $x$ and $y$ are semi-independent these
represent distinct elements of $\Gamma$. In particular we have
$b_m(\{x,y\})\ge2^m$ for each $m>0$, and hence
$\omega(\{x,y\})\ge2$. The assertion now follows from Proposition
\ref{s'n'w prunes}.
\EndProof

\section{Decompositions}\label{decomposition section}

If $T_1$ and $T_2$ are disjoint subtrees of a tree $T$ then there
is a unique open topological arc $A$ which is disjoint from
$T_1$ and $T_2$, and whose endpoints
are vertices, one of which lies in $T_1$ and the other in $T_2$.
We say that an open edge contained in $A$ {\it lies between}
$T_1$ and $T_2$. 

\Proposition\label{XY decomposition}
Let $\Gamma$ be a group and $T$ a $\Gamma$-tree.  Suppose that $x$ and
$y$ are elliptic elements of $\Gamma$ such that $\Fix(x) \cap \Fix(y)
= \emptyset$.  Let $e$ lie between $\Fix(x)$ and $\Fix(y)$.
Suppose that $n$ is a positive integer such that
$e$ is not fixed by $x^k$ nor $y^k$ for $0 < k \le 2n$.
Then there exist disjoint subsets $X$ and $Y$ of $\Gamma$ such that
\begin{itemize}
\item $\Gamma$ is the disjoint union of $X$, $Y$ and $\Gamma_e$;
\item $x^{\pm k}  Y \subset X$ and $y^{\pm k}  X
  \subset Y$ for $0 < k \le n$; and
\item $x^i  Y \cap x^j  Y = \emptyset$ and
$y^iX \cap y^jX = \emptyset$
for any pair of  distinct integers $i$ and $j$
 with $-n\le i \le n$ and $-n\le j\le n$.
\end{itemize}
\EndProposition
\Proof Let $A$ denote the unique open topological arc which is
disjoint from $\Fix(x)\cup\Fix(y)$ and has endpoints $a\in\Fix(x)$ and
$b\in\Fix(y)$.  Thus $e\subset A$.  Let the subtrees $T_x$ and $T_y$
be the components of $T-e$ which contain $a$ and $b$ respectively.
Let $v_x$ and $v_y$ be the endpoints of $e$ which are contained
in $T_x$ and $T_y$ respectively.

We define 
$$\begin{aligned}
 X &= \{\gamma\in\Gamma\;|\;
\gamma\cdot e \subset T_x \};\\
Y &= \{\gamma\in\Gamma\;|\;
\gamma\cdot e \subset T_y \}.
\end{aligned}$$
It is clear that $\Gamma$ is the disjoint union of $X$, $Y$ and
$\Gamma_e$.

Let $k$ be any integer with $0 < k \le n$.  Since $A\supset e$ has
$a\in\Fix(x)$ as an endpoint, and since $x^{\pm k}$ does not fix $e$,
we have $e\not\subset x^{\pm k}\cdot A$. Since $a\in x^{\pm k}\cdot
A$, the component of $T-e$ which contains $x^{\pm k}\cdot A$ is the
one containing $a$, namely $T_x$. This shows that $x^{\pm k}\cdot
A\subset T_x$.  Similarly we have $y^{\pm k}\cdot A\subset T_y$.  In
particular this implies $x^{\pm k}\cdot e\subset T_x$ and
$y^{\pm k}\cdot e\subset T_y$.

Since $T_x$ and $T_y$ are disjoint, we have that
$e\not\subset x^{\pm k}\cdot T_y$ and
$e\not\subset y^{\pm k}\cdot T_x$.  Since the
subtree $x^{\pm k}\cdot T_y$ contains the endpoint $x^{\pm k}\cdot
v_k$ of $x^{\pm k}\cdot e \subset T_x$, it follows from the
connectedness of $T_y$ that $x^{\pm k}\cdot T_y\subset T_x$.
Similarly we have $y^{\pm k}\cdot T_x\subset T_y$.

Now let $\gamma\in Y$ be given.  We have $\gamma\cdot e\subset T_y$
and hence $x^{\pm1}\gamma\cdot e \subset x^{\pm1}\cdot T_y \subset
T_x$.  Hence $x^{\pm1}\gamma\in X$.  This shows that $x^{\pm1} Y
\subset X$.  The proof that $y^{\pm1} X \subset Y$ is similar.

For the proof of the last part of the statement, we assume that $i$
and $j$ are distinct integers with $-n\le i \le n$ and $-n\le j\le n$.
By symmetry it is enough to show that $x^iY$ is disjoint from $x^jY$.
Since $x^k Y = \{\gamma \in \Gamma \;|\; \gamma\cdot e \subset
x^k\cdot T_y\}$, it suffices to show that $x^i\cdot T_y$ is disjoint
from $x^j\cdot T_y$.  Our hypothesis implies that $x^i\cdot e \not=
x^j\cdot e$, and hence that $x^i\cdot v_y$ is distinct from $x^j\cdot
v_y$.  For $k = i, j$, the tree $x^k\cdot T_y$ is the component of
$T-x^k\cdot e$ which contains the vertex $x^k\cdot v_y$.  Since $e$
lies on the geodesic between $a$ and $T_y$, and since $a$ is fixed by
$x$, it follows that the edge $x^j\cdot e$ cannot be contained in
$x^i\cdot T_y$, and the edge $x^i\cdot e$ cannot be contained in
$x^j\cdot T_y$.  In particular $x^i\cdot T_y$ and $x^j\cdot T_y$ are
components of $T - (x^i\cdot e \cup x^j\cdot e)$.  The subtrees
$x^i\cdot T_y$ and $x^j\cdot T_y$ have distinct closest vertices to
$a$, namely $x^i\cdot v_y$ and $x^j\cdot v_y$, so they are not equal.
Therefore they must be disjoint.

\EndProof

The results in \cite{paradoxical} (with Lemma \ref{elles sans foot}
replacing \cite{paradoxical}*{Lemma 5.5} as discussed in Section
\ref{measure section}) will allow us to translate the
combinatorial information provided by Proposition \ref{XY
  decomposition} into measure-theoretic information.

In the following discussion we will regard a Kleinian group $\Gamma$
as acting on the sphere at infinity $S_\infty$ of $\HH^3$. We will use
the same notation as in \cite{paradoxical}.  In particular, the
conformal expansion factor (\cite{paradoxical}*{2.4}) of an element
$\gamma\in\Gamma$ associated to the point $z\in\HH^3$ will be denoted
$\lambda_{\gamma,z}$; the pull-back (\cite{paradoxical}*{3.1}) of a
measure $\mu$ under an isometry $\gamma$ will be denoted
$\gamma^*\mu$; and $\mathcal A = (A_z)$ is the area density on
$S_\infty$ (see \cite{paradoxical}*{3.3}).

\Proposition\label{measure decomposition} Let $M = \HH^3/\Gamma$ be a
closed hyperbolic $3$-manifold and let $P$ be a point in $\HH^3$.
Assume that $T$ is a $\Gamma$-tree and that there are $T$-elliptic
elements $x$ and $y$ of $\Gamma$ such that $\Fix(x) \cap \Fix(y) =
\emptyset$ and neither $x^2$ nor $y^2$ has a fixed edge in $T$.  Then
there exist Borel measures $\nu$, $\nu^+$, $\nu^-$, $\eta$, $\eta^+$,
$\eta^-$ on $S_\infty$ such that
\begin{enumerate}
\item $\nu + \eta \le A_P$;
\item $\nu^+ + \nu^- \le \nu$ and  $\eta^+ + \eta^- \le \eta$  
\item $\int_{S_\infty}\lambda^2_{x^{\pm1}, P}d\nu^\pm =
  \eta(S_\infty)$ and $\int_{S_\infty}\lambda^2_{y^{\pm1}, P}d\eta^\pm =
  \nu(S_\infty)$ 
\end{enumerate}
Furthermore, in the case that $T$ is the dual tree of an
incompressible surface $F$ in the closed $3$-manifold $\HH^3/\Gamma$,
we have $\nu + \eta = A_P$.
\EndProposition

\Proof Note that the hypotheses of Proposition
\ref{XY decomposition} hold with $n=1$. Let $X$ and $Y$ be the subsets of $\Gamma$ given by Proposition
\ref{XY decomposition}.  Thus $\Gamma$ is the disjoint union of $X$,
$Y$ and $\Gamma_e$ for some edge $e$ of $T$. Furthermore, $x^{-1}Y$ and $xY$ are mutually disjoiut subsets of $X$, while $y^{-1}X$ and $yX$ are mutually disjoiut subsets of $X$. Since $\Gamma$ is
discrete, the set $\Gamma\cdot P$ is uniformly discrete, in the sense
of \cite{paradoxical}*{4.1}.
Set $X' = X - (xY\cup x^{-1}Y)$ and $Y' = Y- (yX\cup y^{-1}X)$.
Define $\mathfrak V$ to be the collection of
all unions of sets in
$$\{X'\cdot P,\; Y'\cdot P,\;
xY\cdot P,\; x^{-1}Y\cdot P,\;
yX\cdot P,\; y^{-1}X\cdot P,\;
 \Gamma_e\cdot P
\}.$$ 

We apply \cite{paradoxical}*{Proposition 4.2} with $W = \Gamma\cdot
P$, to construct a family $(\calm_{V})_{V\in \mathfrak V}$, of
$D$-dimensional conformal densities, for some $D\in [0,2]$, such that
conditions (i)-(iv) of \cite{paradoxical}*{Proposition 4.2} hold.  Set
$\calm_{X\cdot P} = (\nu_z)$, $\calm_{xY\cdot P} = (\nu^+_z)$,
$\calm_{x^{-1}Y\cdot P} = (\nu^-_z)$, $\calm_{Y\cdot P} = (\eta_z)$,
$\calm_{yX\cdot P} = (\eta^+_z)$, $\calm_{y^{-1}X\cdot P} =
(\eta^-_z)$, and $\calm_{\Gamma_e\cdot P} = (\epsilon_z)$.  It follows
from conditions (i) and (ii) of \cite{paradoxical}*{Proposition 4.2}
that $\calm_\Gamma = \calm_{X\cdot P} + \calm_{Y\cdot P} + \mathcal
M_{\Gamma_e\cdot P}$ is a $\Gamma$-invariant conformal density.  Since
$M = \HH^3/\Gamma$ is a closed manifold, every $\Gamma$-invariant
superharmonic function on $M$ is constant.  Thus by
\cite{paradoxical}*{Proposition 3.9}, $D=2$ and $\calm_{X\cdot P} +
\mathcal M_{Y\cdot P} + \calm_{\Gamma_e\cdot P} = k\mathcal A$ for
some constant $k$.  Condition (i) of \cite{paradoxical}*{Proposition
  4.2} guarantees that $k > 0$.  Thus by normalizing the family
$(\mathcal M_{V})_{V\in \mathfrak V}$ appropriately we may assume that
$k=1$.

We define $\nu = \nu_P$, $\nu^\pm=\nu^\pm_P$, $\eta = \eta_P$, and
$\eta^\pm=\eta^\pm_P$.  Then Conclusion (1) follows from the equality
$\calm_{X\cdot P} + \mathcal M_{Y\cdot P} + \calm_{\Gamma_e\cdot P} =
k\mathcal A$ by specializing to $z=P$.

According to \cite{paradoxical}*{Proposition 4.2 (ii)}, we have 
$$\calm_{X\cdot P} = \calm_{X'\cdot P} + \calm_{xY\cdot P}
+M_{x^{-1}Y\cdot P}$$
and 
$$\calm_{Y\cdot P} = \calm_{Y'\cdot P} +
\calm_{yX\cdot P} + \calm_{y^{-1}X\cdot P}.$$
Specializing to $z=P$, we obtain
Conclusion (2).

Now we turn to the proof of Conclusion (3).  We will only give the
proof that $\int_{S_\infty}\lambda^2_{x,P}d\nu^+ = \eta(S_\infty)$;
the other three statements included in Conclusion (3) are proved by
the same argument.  The proof is based on
\cite{paradoxical}*{Proposition 4.2 (iii)}, which asserts that
$x^*_\infty(\calm_{xY}) = \calm_Y$.  After substituting
$\calm_{xY\cdot P} = (\nu^+_z)$ and $\calm_{Y\cdot P} = (\eta_z)$, the
definition of the pull-back (see \cite{paradoxical}*{3.4.1}) states
that $d(x^*_\infty\nu^+_z) = \lambda^2_{x,z}d\nu^+_z$.  Taking $z=P$
and integrating over $S_\infty$ we obtain
$$\int_{S_\infty}\lambda^2_{x,P}d\nu^+ =
 \int_{S_\infty}d(x^*_\infty\nu^+) =
 \int_{S_\infty}d\eta = \eta(S^\infty)$$
as required for Conclusion (3).

We now consider the case where $T$ is the dual tree of an
incompressible surface $F$ in the $3$-manifold $\HH^3/\Gamma$.  For $z
\in \HH^3$ the support of the measure $\epsilon_z$ is contained in the
limit set of the Kleinian group $\Gamma_e$.  The group $\Gamma_e$
cannot be conjugate to the fundamental group of a fiber or a
semi-fiber in $M$, since in that case the tree $T$ would be a line and
the square of any elliptic element of $\Gamma$ would fix the entire
tree $T$.  It follows from work of Thurston, Bonahon and Canary (see
\cite{canary}*{Corollary 8.1}) that the group $\Gamma_e$ is either
geometrically finite or a virtual fiber subgroup.  If $\Gamma_e$ were
a virtual fiber subgroup then, since $\Gamma_e$ is equal to the image
of $\pi_1(F)$ up to conjugacy, we would deduce, by applying
\cite{hempel}*{Theorem 10.5} to the components of the manifold
obtained by splitting $M$ along $F$, that $F$ is a fiber or
semi-fiber.  Thus $\Gamma_e$ must be geometrically finite, and by a
theorem of Ahlfors's \cite{ahlfors} the limit set of $\Gamma_e$ has
area measure 0.  Since the measure $\epsilon_P$ is supported on the
limit set of $\Gamma_e$, we see that $\epsilon_P$ is singular with
respect to $A_P$ and we have $\nu + \eta = \nu_P + \eta_P = A_P$, as
required for the last sentence of the statement.
\EndProof

\Proposition\label{double trouble} Let $M = \HH^3/\Gamma$ be a closed
hyperbolic $3$-manifold and let $P$ be a point of $\HH^3$.  Suppose
that $T$ is the dual $\Gamma$-tree of an incompressible surface in
$M$, and that $x$ and $y$ are $T$-elliptic elements of $\Gamma$ such
that $\Fix(x) \cap \Fix(y) = \emptyset$ and neither $x^2$ nor $y^2$
has a fixed edge in $T$.  Set $D_x=\exp(2d_P(x))$ and
$D_y=\exp(2d_P( y))$. Suppose that $\nu$, $\eta$, $\nu^+$,
$\nu^-$, $\eta^+$, $\eta^-$ are measures having the properties stated
in (the conclusion of) Proposition \ref{measure decomposition}, and
let $\alpha$, $\beta$, $\alpha^+$, $\alpha^-$, $\beta^+$ and $\beta^-$
denote their respective total masses.  Then we have
\begin{enumerate}
\item$\alpha+\beta=1$;
\item$\displaystyle\frac{\beta (1-\alpha^+)}{\alpha^+(1-\beta )} \le D_x$\;; \quad
$\displaystyle\frac{\beta (1-\alpha^-)}{\alpha^-(1-\beta )} \le D_x$\;; 
\vspace{6pt plus 3pt}
\item$\displaystyle\frac{\alpha (1-\beta ^+)}{\beta^+(1-\alpha )} \le D_y$\;; \quad
$\displaystyle\frac{\alpha (1-\beta^-)}{\beta^-(1-\alpha )} \le D_y$.
\end{enumerate} 
\EndProposition
 
\Proof Since $T$ is the dual $\Gamma$-tree of an
incompressible surface, the last sentence of Proposition \ref{measure
  decomposition} gives that $\nu+\eta=A_P$, which implies (1). We
intend to deduce the inequalities in (2) and (3) from Lemma
\ref{elles sans foot}. For (2), the element $\gamma$ in the
statement of Lemma \ref{elles sans foot} should be replaced by
$x^{\pm1}$, the measure $\nu$ by $\nu^{\pm}$, $a$ by $\alpha^{\pm}$
and $b$ by $\beta$.  For (3), we should replace $\gamma$ by
$y^{\pm1}$, $\nu$ by $\eta$, $a$ by $\beta^\pm$ and $b$ by $\alpha$.
The only obstruction to this argument is that to apply Lemma
\ref{elles sans foot} we must ensure that we do not have $a=b=0$ or
$a=b=1$ in any of these four cases.

We claim that $\alpha \not= 0$.  Otherwise we would have $\alpha^+=0$,
so the measure $\nu^+$ would be singular.  But by Proposition
\ref{measure decomposition} we then would have
$$\beta = \eta(S_\infty) = \int_{S_\infty}\lambda^2_{x, P}d\nu^+ = 0,$$
contradicting the fact that $\alpha + \beta = 1$.  Similarly, $\beta\not=0$.
In addition, since $\alpha + \beta = 1$, we cannot
have $\alpha^\pm = \beta = 1$ nor $\alpha = \beta^\pm = 1$.
This shows that the possibilities $a=b=0$ or $a=b=1$ do not
arise, and the argument given above does, in fact, prove the Proposition.
\EndProof

\Proposition\label{long squares} Let $M = \HH^3/\Gamma$ be a closed
hyperbolic $3$-manifold and let $P$ be a point of $\HH^3$.  Suppose
that $T$ is the dual $\Gamma$-tree of an incompressible surface in
$M$, and that $x$ and $y$ are $T$-elliptic elements of $\Gamma$ such
that $\Fix(x) \cap \Fix(y) = \emptyset$ and neither $x^2$ nor $y^2$
has a fixed edge in $T$.  Set $D_x=\exp(2d_P(x))$ and
$D_y=\exp(2d_P(y))$. Then
$$\frac{\Xsqrt{8D_x+1}-3}{D_x-1} + \frac{\Xsqrt{8D_y+1}-3}{D_y-1} \le 2.$$
\EndProposition
\Proof
Let $\alpha$, $\beta $, $\alpha^+$, $\alpha^-$, $\beta^+$, $\beta^-$ in $[0,1]$
be the numbers defined as in Proposition \ref{double trouble}.
By symmetry may assume that
$\alpha^+\le \alpha/2$ and $\beta^+\le \beta /2$.  Since the function
$f(t) = (1-t)/t$ is decreasing on the interval $(0,1)$, it follows from Proposition \ref{double trouble} that
$$
D_x \ge \left(\frac{\beta }{1-\beta }\right)\left(\frac{1-\alpha/2}{\alpha/2}\right) =
\left(\frac{\beta }{1-\beta }\right)\left(\frac{2-\alpha}{\alpha}\right) = 
\frac{(1-\alpha)(2-\alpha)}{\alpha^2} 
$$
and
$$
D_y \ge \left(\frac{\alpha}{1-\alpha}\right)\left(\frac{1-\beta /2}{\beta /2}\right) =
\left(\frac{\alpha}{1-\alpha}\right)\left(\frac{2-\beta }{\beta }\right) =
\frac{(1-\beta )(2-\beta )}{\beta ^2} .
$$

Let us set $\overline D_x = (1-\alpha)(2-\alpha)/\alpha^2$ and
$\overline D_y = (1-\beta )(2-\beta )/\beta ^2$.  Solving the
quadratic equations, we find that
$$
2\alpha = \frac{\Xsqrt{8\overline D_x + 1} -3}{\overline D_x - 1} ;
\text{ and }
2\beta  = \frac{\Xsqrt{8\overline D_y + 1} -3}{\overline D_y - 1} .
$$
A straightforward computation shows that the function $g(t) =
(\sqrt{8t + 1} -3)/(t - 1)$ is decreasing on the interval $(1,
\infty)$.  Since $\overline D_x \le D_x$ and $\overline D_y \le D_y$
we thus have
$$\frac{\Xsqrt{8D_x+1}-3}{D_x-1} + \frac{\Xsqrt{8D_y+1}-3}{D_y-1} \le 2(\alpha+\beta ) = 2.$$  
\EndProof

\Corollary Let $M = \HH^3/\Gamma$ be a closed hyperbolic $3$-manifold
and let $P$ be a point of $\HH^3$.  Suppose that $T$ is the dual
$\Gamma$-tree of an incompressible surface in $M$, and that $x$ and
$y$ are $T$-elliptic elements of $\Gamma$ such that $\Fix(x) \cap
\Fix(y) = \emptyset$ and neither $x^2$ nor $y^2$ has a fixed edge in
$T$.  Then
$$\max(d_P(x), d_P(y)) \ge \frac12\log 3.$$
\EndCorollary

\Proof
According to Proposition \ref{long squares} we have
$$\frac{\Xsqrt{8D_x+1}-3}{D_x-1} + \frac{\Xsqrt{8D_y+1}-3}{D_y-1} \le 2,$$ 
where $D_x=\exp(2d_P(x))$ and
$D_y=\exp(2d_P(y))$.  In view of the symmetry of the desired conclusion, it suffices to consider the case in which
$$\frac{\Xsqrt{8D_x+1}-3}{D_x-1} \le 1.$$
Solving, we find $D_x \le 3$, and hence $d_P(x) \ge \frac12\log 3$, from which the conclusion follows in this case.
\EndProof

\section{Some hyperbolic trigonometry}\label{sphere section}

We shall denote by $d_s$ the spherical distance on the unit sphere
$S^2\subset\EE^3$.

\Proposition\label{mademoiselle victoire}
If $\eta_1,\ldots,\eta_n$ are points on $S^2$, we have
$$\sum_{1\le i<j\le n}\cos d_s(\eta_i,\eta_j)\ge-n/2.$$
\EndProposition

\Proof We regard $S^2$ as the unit sphere in $\EE^3$, and we let
$v_i\in\RR^3$ denote the position vector of $\eta_i$. We have
$$\begin{aligned}
0&\le\langle\sum_{i=1}^nv_i,\sum_{i=1}^nv_i\rangle\\
 &=\sum_{i=1}^n\|v_i\|^2+\sum_{i\ne j}\langle v_i,v_j\rangle\\
 &=n+2\sum_{1\le i< j\le n}\langle v_i,v_j\rangle\\
 &=n+2\sum_{1\le i< j\le n}\cos d_s(\eta_i,\eta_j),
\end{aligned}$$
from which the conclusion follows.
\EndProof

\Corollary\label{scholastique}Let $P$ be a point in $\HH^3$, and let $Q_1,\ldots,Q_n\in\HH^3$ be points distinct from $P$. Then we have
$$\sum_{1\le i<j\le n}\cos \angle(Q_i,P,Q_j)\ge-n/2.$$
\EndCorollary

\Proof
We consider 
the unit sphere $\Sigma$  in the tangent space to $\HH^3$ at $P$. For $i=1,\ldots,n$ we let $r_i$ denote the ray from $P$ to $Q_i$, and let
$\eta_i\in \Sigma $ denote the unit tangent vector to the ray $r_i$. Then for
any two distinct indices $i,j\in\{1,\ldots,n\}$, we have
$d_s(\eta_i,\eta_j)=\angle (Q_i,P,Q_j).$ The conclusion now follows from Proposition \ref{mademoiselle victoire}.
\EndProof

\Proposition\label{alhambra}Let $\nu$ be a positive real number, and
let $Q,R,S$ be points of $\HH^n$, for some $n\ge2$, such that
$\max(\dist(Q,R),\dist(R,S))\le\nu$. Then
$$ \dist(Q,S)\le\max(\nu,\arccosh(\cosh^2\nu-\cos(\angle QRS)\sinh^2\nu)).$$
\EndProposition

\Proof
Let $q$ and $s$ denote the rays that originate at $R$ and pass,
respectively, through $Q$ and $S$. Let $Q'$ and $S'$ denote the points
that lie on $q$ an $r$ respectively, and have distance $\nu$ from
$R$. Then $\dist(Q,S)$ is bounded above by the diameter of the
triangle $\Delta$ with vertices $Q'$, $R$ and $S'$. This
diameter is in turn equal to the maximum of the side lengths of
$\Delta$. Two of these side lengths are equal to $\nu$, and the third
is equal to $\arccosh(\cosh^2\nu-\cos(\angle
QRS)\sinh^2\nu))$ by the hyperbolic law of cosines.
\EndProof

\Lemma\label{d'alembert}Let $x$ and $y$ be isometries of
$\HH^3$. Let $P$ be a point of $\HH^3$, and let $\nu$ be a positive
real number. Suppose that
$$d_P(x) \le \nu< d_P(x^2)$$
 and that
$$d_P(y) \le \nu< d_P(y^2).$$
Set 
$$A=\frac{\cosh^2\nu-\cosh\nu}{\sinh^2\nu}.$$
Then $$\sum_{(u,v)\in \{\pm1\}\times \{\pm1\}}\cos\angle(x^u\cdot P,P,y^v\cdot P)>-2-2A.$$
\EndLemma

\Proof
Noting that $d_P(x^{-1}) =d_P(x) \le \nu$, and applying Proposition
\ref{alhambra} with $Q=x^{-1}\cdot P$, $R=P$ and $S=x\cdot P$, we find
that
$$\dist(x^{-1}\cdot P,x\cdot P)\le \max(\nu,
\arccosh (
\cosh^2\nu-\cos(\angle (x^{-1}\cdot P,P,x\cdot P))\sinh^2\nu)).$$
On the other hand, using the hypothesis we find that
$$\dist(x^{-1}\cdot P,x \cdot P)=\dist(x^2 \cdot P,P)>\nu.$$
Hence
$\nu<
\arccosh (
\cosh^2\nu-\cos(\angle (x^{-1}\cdot P,P,x\cdot P))\sinh^2\nu)$, i.e.
$$\cosh\nu< \cosh^2\nu-\cos(\angle (x^{-1}\cdot P,P,x\cdot P))\sinh^2\nu.$$
In view of the definition of $A$, this implies that
\Equation\label{eggs case}
\cos\angle (x^{-1}\cdot P,P,x\cdot P)< A.
\EndEquation
The same argument shows that
\Equation\label{why case}
\cos\angle (y^{-1}\cdot P,P,y\cdot P)< A.
\EndEquation

Now set $\xi_1=x$, $\xi_2=x^{-1}$, $\xi_3=y$ and $\xi_4=y^{-1}$. It
follows from Corollary \ref{scholastique}, with $Q_i=\xi_i\cdot P$,
that \Equation\label{UMbrella} \sum_{1\le i<j\le n}\cos
\angle(\xi_i\cdot P,P,\xi_j\cdot P )\ge-n/2.
\EndEquation
Note that the left-hand sides of (\ref{eggs case}) and (\ref{why
  case}) are among the six terms on the left side of
(\ref{UMbrella}). The remaining terms are the quantities
$\cos\angle(x^u\cdot P,P,y^v\cdot P)$, where $(u,v)$ ranges over
$\{\pm1\}\times \{\pm1\}$. Hence the conclusion of the lemma follows
from (\ref{eggs case}), (\ref{why case}) and (\ref{UMbrella}).
\EndProof

\Notation\label{phidef} We define a function $\phi$
on $(0,\infty)$ by
$$\phi(t)=\max(3t, 2\arccosh \big(2\cosh^2t-\frac12\cosh t-\frac12\Big)).  
$$
Note that the expression $\arccosh
\big(2\cosh^2\nu-\dfrac12\cosh\nu-\dfrac12\big)$ is well-defined and
positive for $\nu>0$, since we have $2t^2-\dfrac12t-\dfrac12>1$ for
$t>1$. Furthermore, since $2t^2-\dfrac12t-\dfrac12>1$ is strictly
monotone increasing for $t>1$, the function $\phi$ is also strictly
monotone increasing on $(0,\infty)$.
\EndNotation

\Lemma\label{mlle. hus}Let $x$ and $y$ be isometries of
$\HH^3$. Let $P$ be a point of $\HH^3$, and let $\nu$ be a positive
real number. Suppose that
$$d_P(x) \le \nu < d_P(x^2)$$
 and that
$$d_P(y) \le \nu < d_P(y^2).$$
Then
$$\min(d_P(xy)+d_P(yx),d_P(xy^{-1})+d_P(y^{-1}x))
\le \phi(\nu).$$
\EndLemma

\Proof We set $$E= d_P(xy)+d_P(yx)=\dist(x^{-1}\cdot P ,y\cdot
P)+\dist(x\cdot P ,y^{-1}\cdot P)$$
and $$E'=d_P(xy^{-1})+d_P(y^{-1}x)=\dist(x^{-1}\cdot P ,y^{-1}\cdot
P)+\dist(x\cdot P,y\cdot P).$$
We are required to prove that 
\Equation\label{muriel}
\min(E,E')
\le \phi(\nu).
\EndEquation

Since $\max(d_P(x),d_P(y)) <\nu$, each of the terms $d_P(xy)$,
$d_P(yx)$, $d_P(xy^{-1})$ and $d_P(y^{-1}x)$ is bounded above by
$2\nu$. If one of these terms is bounded above by $\nu$, then the left
hand side of (\ref{muriel}) is bounded above by
$3\nu\le\phi(\nu)$. Hence we may assume that each of these four terms
is greater than $\nu$.

Set 
$$A=\frac{\cosh^2\nu-\cosh\nu}{\sinh^2\nu}.$$
According to Lemma \ref{d'alembert} we have
\Equation\label{same old same old}
\sum_{(u,v)\in \{\pm1\}\times \{\pm1\}}\cos\angle(x^u\cdot P,P,y^v\cdot P)>-2-2A.
\EndEquation
The left hand side of (\ref{same old same old}) may be written as $C+C'$, where
$$C=\cos\angle(x\cdot P,P,y^{-1}\cdot P) +\cos\angle(x^{-1}\cdot P,P,y\cdot P) $$
and
$$C'=\cos\angle(x\cdot P,P,y\cdot P) +\cos\angle(x^{-1}\cdot P,P,y^{-1}\cdot P). $$ 
In particular we have
\Equation\label{Sandy's pants}
\max(C,C')\ge-1-A.
\EndEquation

Consider the case in which $C\ge-1-A$. Applying Proposition
\ref{alhambra} with $Q=x\cdot P$, $R=P$ and $S=y^{-1}\cdot P$, we find
that
$$
d(x\cdot P,y^{-1}\cdot P)\le
\max(\nu,
\arccosh(
\cosh^2\nu-\cos(\angle (x\cdot P,P,y^{-1}\cdot P))\sinh^2\nu).
$$
Since $d(x\cdot P,y^{-1}\cdot P)>\nu$ it follows that
\Equation\label{pee-wee reese} \cosh d(x\cdot P,y^{-1}\cdot P)\le
\cosh^2\nu-\cos(\angle (x\cdot P,P,y^{-1}\cdot P))\sinh^2\nu.
\EndEquation
Similarly,
\Equation\label{junior gilliam}
\cosh d(x^{-1}\cdot P,y\cdot P)
\le 
\cosh^2\nu-\cos(\angle (x^{-1}\cdot P,P,y\cdot P))\sinh^2\nu).
\EndEquation
Adding (\ref{pee-wee reese}) and (\ref{junior gilliam}), and using the definition of $C$, we obtain
\Equation\label{duke snider}
\begin{aligned}\cosh d(x\cdot P,y^{-1}\cdot P)+\cosh d(x^{-1}\cdot P,y\cdot P)
&\le
2\cosh^2\nu-C\sinh^2\nu\\
&\le
2\cosh^2\nu+(1+A)\sinh^2\nu\\
&=4\cosh^2\nu-\cosh\nu-1.\end{aligned}
\EndEquation
On the other hand, since $\cosh$ is convex, we have

$$
\begin{aligned}
\cosh (E/2)&=\cosh(\frac12(d(x^{-1}\cdot P ,y\cdot P)+d(x\cdot P ,y^{-1}\cdot P)))\\
&\le\frac12(\cosh d(x\cdot P,y^{-1}\cdot P)+\cosh d(x^{-1}\cdot P,y\cdot P)),
\end{aligned}
$$
which with (\ref{duke snider}) gives 
$$\cosh (E/2)\le 2\cosh^2\nu-\frac12\cosh\nu-\frac12.$$
This implies (\ref{muriel}).

If $C'\ge-1-A$, the same argument shows that
$$\cosh (E'/2)\le  2\cosh^2\nu-\frac12\cosh\nu-\frac12,$$
which again implies (\ref{muriel}).

Thus, in view of (\ref{Sandy's pants}), the conclusion is seen to hold
in all cases.
\EndProof

\section{Proof of the main theorem}\label{putting it together}

\Lemma\label{commute} Let $\Gamma$ be the fundamental group of an
orientable hyperbolic $3$-manifold. Then the centralizer of every
non-trivial element of $\Gamma$ is abelian.  Furthermore, if $t$ and
$u$ are elements of $\Gamma$ and if $t$ commutes with $utu^{-1}$, then
$t$ commutes with $u$.
\EndLemma

\Proof Up to isomorphism, we may identify $\Gamma$ with a torsion-free
discrete subgroup of $\isomplus(\HH^3)$. Any non-trivial element $x$
of $\Gamma$ is either loxodromic or parabolic. In these respective
cases we let $A_x$ denote the axis of $x$ or its fixed point on the
sphere at infinity.

If $1\ne x\in\Gamma$, any element of the centralizer of $x$ must leave
$A_x$ invariant. Since $\Gamma$ is discrete and torsion-free, the
stabilizer of $A_x$ in $\Gamma$ is abelian. This proves the first
assertion.

In proving the second assertion we may assume that $t$ and $u$ are
non-trivial. We have $u\cdot A_t=A_{utu^{-1}}$. On the other hand,
since $utu^{-1}$ commutes with $t$, we have $A_{utu^{-1}}=A_t$. Hence
$u\cdot A_t=A_t$, so that $u$ leaves $A_t$ invariant and therefore
commutes with $t$.
\EndProof

\Corollary\label{commutators}  Let $\Gamma$ be the fundamental group of an
 orientable hyperbolic $3$-manifold,
and let $x$ and $y$ be elements of $\Gamma$.
\begin{enumerate}
\item\label{langevin} If $x$ and $y$ do not commute, then
  $xyx^{-1}y^{-1}$ and $yx^{-1}y^{-1}x$ do not commute.
\item\label{lauvergnat} If $x^m$ and $y^n$ commute, for some non-zero
  integers $m$ and $n$, then $x$ and $y$ commute.
\end{enumerate}
\EndCorollary

\Proof To prove (\ref{langevin}), suppose that $xyx^{-1}y^{-1}$ commutes with $yx^{-1}y^{-1}x$.
Apply Lemma \ref{commute} with $t = xyx^{-1}y^{-1}$ and $u = x^{-1}$
to deduce that $yx^{-1}y^{-1}$ commutes with $x^{-1}$.  Then apply the
lemma again, with $t = yx^{-1}y^{-1}$ and $u = y$, to deduce that $x$
commutes with $y$, a contradiction.

In proving (\ref{lauvergnat}) , we may assume that $x\ne1$ and
$y\ne1$. Since $\Gamma$ is torsion-free it then follows from Lemma
\ref{commute} that for any $m\ne0$ the centralizer $C$ of $x^m$ is
abelian. We have $x\in C$, and if $x^m$ commutes with $y^n$ for some
$n\ne0$ then $y^n\in C$.  Since $C$ is abelian it follows that $x$
commutes with $y^n$.  But Lemma \ref{commute} also implies that the
centralizer $C'$ of $y^n$ is abelian, and since $C'$ contains $x$ and
$y$ we conclude that $x$ and $y$ commute, as required.
\EndProof

\Proposition\label{sanity clause} Let $M = \HH^3/\Gamma$ be a closed
hyperbolic $3$-manifold containing an incompressible surface $F$ and
let $T$ denote the dual $\Gamma$-tree of $F$. Let $x$ and $y$ be
non-commuting elements of $\Gamma$, and let $P$ be a point of $\HH^3$.
\begin{enumerate}
\item \label{won}If $x$ and $y$ are both $T$-hyperbolic, and $F$ is not a fiber
  or a semi-fiber, then
$$\max(d_P(x),d_P(y))\ge\frac12\log2=0.346\ldots.$$
\item\label{too} If $x^2$ fixes at least one edge of $T$, then
 $$\max(d_P(x),d_P(y))\ge\frac12\log\alpha=0.304\ldots,$$
where $\alpha=1.839\ldots$ is the unique real root of the polynomial
$Q(t)=t^3-t^2-t-3$.
\item\label{fore}If $x$ is $T$-elliptic and $\Fix(x)\cap\Fix(yxy^{-1})\ne\emptyset$, then 
$$\max(d_P(x),d_P(y))\ge\log\gamma=0.593\ldots,$$
where $\gamma=1.8105\ldots$ is the unique real root of the polynomial
$R(t)=t^4-t^3-t-3$.
\item\label{tree} If $x$ is $T$-elliptic 
then 
$$\max(d_P(x),d_P(y))\ge0.286.$$
\end{enumerate}
\EndProposition

\Proof To prove (\ref{won}) we observe that if $F$ is not a fiber or a
semi-fiber, then the action of $\Gamma$ on $T$ is non-trivial by
Proposition \ref{nontrivial tree} and linewise faithful by Proposition
\ref{yes georgia}. Hence by Proposition \ref{ILGWU}, if $x$ and $y$
are both $T$-hyperbolic then they are semi-independent. The assertion
therefore follows from Corollary \ref{florida}.

To prove (\ref{too}), we first observe that since $x$ and $y$ do not
commute, it follows from Assertion (\ref{lauvergnat}) of Corollary
\ref{commutators} that $x^2$ and $y$ do not commute.  It now follows
from Lemma \ref{commute} that $x^2$ and $yx^2y^{-1}$ do not commute.
If $x^2$ fixes at least one edge of $T$, then $yx^2y^{-1}$ also fixes
at least one edge of $T$, and hence by Proposition \ref{arizona},
$x^2$ and $yx^2y^{-1}$ are independent in $\Gamma$. It therefore
follows from \cite{acs}*{Theorem 4.1} (which is in turn a consequence
of results proved in \cite{accs}, \cite{agol} and \cite{cg}) that
$$\frac1{1+\exp d_P(x^2)}+\frac1{1+\exp d_P(yx^2y^{-1})}\le\frac12.$$
If we set $D=\max(d_P(x),d_P(y))$, we have 
$d_P(x^2)\le2D$ and
$d_P(yx^2y^{-1})\le4D$. Hence 
$$\frac1{1+\exp(2D)}+\frac1{1+\exp(4D)}\le\frac12.$$
If we now set $u=\exp (2D)$ we obtain $u^3-u^2-u-3\ge0$. But the
polynomial $Q(t)=t^3-t^2-t-3$ increases monotonically for
$t\ge1$. Hence $\exp(2D)=u\ge\alpha$, and the conclusion follows.

To prove (\ref{fore}) we note that, by Proposition
\ref{nontrivial tree} and Proposition \ref{old mexico},
the elements $x$ and $yxy^{-1}$ of $\Gamma$ are independent. Hence by
\cite{acs}*{Theorem 4.1} we have
$$\frac1{1+\exp d_P(x)}+\frac1{1+\exp d_P(yxy^{-1})}\le\frac12.$$
If we set $D=\max(d_P(x),d_P(y))$, we have 
$d_P(x)\le2D$ and
$d_P(yxy^{-1})\le3D$. Hence 
$$\frac1{1+\exp D}+\frac1{1+\exp(3D)}\le\frac12.$$
If we now set $v=\exp D$ we obtain $v^4-v^3-v-3 \ge0$. But the
polynomial $R(t)=t^4-t^3-t-3$ increases monotonically for
$t\ge1$. Hence $\exp D=v\ge\gamma$, and the conclusion follows.

To prove (\ref{tree}) we first consider the special case in which the
inequality \Equation\label{squares have more fun}
\max(d_P(x),d_P(y))<\min(d_P(x^2),d_P(y^2))
\EndEquation
holds. 

In the subcase where $x^2$ fixes at least one edge of $T$, the
assertion follows from assertion (\ref{too}), which we have already
proved. Likewise, in the subcase where
$\Fix(x)\cap\Fix(yxy^{-1})\ne\emptyset$, the assertion follows from
assertion (\ref{fore}). We therefore need only address the subcase in
which $\Fix(x)\cap\Fix(yxy^{-1})=\emptyset$ and $x^2$ fixes no edge of
$T$.  Note that in this subcase, the hypotheses of Proposition
\ref{long squares} hold with the elements $x$ and $yxy^{-1}$ of
$\Gamma$ playing the respective roles of $x$ and $y$ in the latter
proposition.

In this subcase we shall assume that $\max(d_P(x),d_P(y))<0.286$ and
deduce a contradiction. In view of (\ref{squares have more fun}), we
may choose a real number $\nu<0.286$ such that
$$\max(d_P(x),d_P(y))<\nu<\min(d_P(x^2),d_P(y^2)).$$
In particular the hypotheses of Lemma \ref{mlle. hus} hold with this
choice of $\nu$, and it follows from that lemma that
$\min(d_P(xy)+d_P(yx),d_P(xy^{-1})+d_P(y^{-1}x)) \le \phi(\nu)$, where
$\phi$ is the function defined in \ref{phidef}.  In particular we have
$$\min(d_P(xy),d_P(yx),d_P(xy^{-1}),d_P(y^{-1}x))
\le \frac12\phi(\nu).$$
If $d_P(xy) \le \dfrac12\phi(\nu)$ then
$$d_P(y^{-1}xy)\le d_P(y^{-1})+d_P(xy)=d_P(y)+d_P(xy) \le \dfrac12\phi(\nu)+\nu.$$
Similarly, if $d_P(xy^{-1}) \le \dfrac12\phi(\nu)$ then $d_P(y^{-1}xy)\le \dfrac12\phi(\nu)+\nu$, and if
$d_P(yx) \le \dfrac12\phi(\nu)$ or $d_P(y^{-1}x)) \le \dfrac12\phi(\nu)$ then
$d_P(yxy^{-1})\le \dfrac12\phi(\nu)+\nu$. Hence,
after possibly interchanging the roles of $y$ and $y^{-1}$, we may assume that
$$d_P(yxy^{-1}) \le \dfrac12\phi(\nu)+\nu.$$

In view of the monotonicity of $\phi$ (see \ref{phidef}), we have
\Equation\label{hang your heart}
d_P(yxy^{-1})\le\dfrac12\phi(0.286)+0.286<0.8227.
\EndEquation

If we set $D_x=\exp(2d_P(x))$ and $D_{yxy^{-1}}=\exp(2d_P(
yxy^{-1}))$, we have $D_x\le\exp(2\cdot0.286)<1.772$ and
$D_y\le\exp(2\cdot0.8227)<5.1831$. Applying Proposition \ref{long
  squares} with the elements $x$ and $yxy^{-1}$ of $\Gamma$ playing
the respective roles of $x$ and $y$ in the latter proposition, and noting
that the function $g(t) = (\sqrt{8t + 1} -3)/(t - 1)$ is strictly
decreasing on the interval $(1, \infty)$, we find that
$$\begin{aligned}2&\ge g(D_x)+g(D_{yxy^{-1}})\\
&>g(1.772)+g(5.1831)\\
&=2.0007\ldots,
\end{aligned}$$
which is the required contradiction in this case.

We now turn to the general case of (\ref{tree}), in which the
inequality \ref{squares have more fun} is not assumed to hold. Again
we argue by contradiction, assuming that $\max(d_P(x),d_P(y))<0.286.$
Since $x$ and $y$ are loxodromic, the quantities $d_P(x^n)$ and
$d_P(y^n)$ tend to $\infty$ with $n$. Hence there is a largest integer
$n_1$ such that $d_P(x^{n_1})<0.286$, and there is a largest integer
$n_2$ such that $d_P(y^{n_2})<0.286$. If we set $x_0=x^{n_1}$ and
$y_0=y^{n_2}$, it follows that $d_P(x_0^2)\ge0.286$ and that
$d_P(y_0^2 )\ge0.286$. Hence
$$\min(d_P(x_0^2),d_P(y_0^2))\ge0.286>
\max(d_P(x_0),d_P(y_0)).$$ The element $x_0=x^{n_1}$ of $\Gamma$ is
$T$-elliptic since $x$ is $T$-elliptic, and since $x$ and $y$ do not
commute it follows from Assertion (\ref{lauvergnat}) of Corollary
\ref{commutators} that $x_0=x^{n_1}$ and $y_0=y^{n_2}$ do not
commute. By the special case of (\ref{tree}) already proved, with
$x_0$ and $y_0$ playing the roles of $x$ and $y$, it now follows that
$\max(d_P(x_0),d_P(y_0))\ge0.286$. This is a contradiction.
\EndProof

\Remark Conclusions (\ref{too}) and (\ref{fore}) of Proposition
\ref{sanity clause} could be improved by using Lemma \ref{mlle. hus},
but this would not affect our main result in this paper.
\EndRemark

\Proposition\label{off the wall} Let $M$ be a hyperbolic $3$-manifold
such that {\it either} $H_1(M;\QQ)\ne0$ {\it or} $M$ is closed and
contains a semi-fiber. Then $0.292$ is a Margulis number for $M$.
\EndProposition

\Proof We write $M=\HH^3/\Gamma$. We suppose that $x$ and $y$ are
elements of $\Gamma$ and that $P$ is a point of $\HH^3$ such that
$\max(d_P(x),d_P(y))<0.292$. We must show that $x$ and $y$ commute.

We first consider the special case in which the inequality
\Equation\label{squares still have more fun}
\min(d_P(x^2),d_P(y^2))\ge0.292
\EndEquation
holds.  In this case the hypotheses of Lemma \ref{mlle. hus} hold with
$\nu=0.292$, and it follows from that lemma that
$\min(d_P(xy)+d_P(yx),d_P(xy^{-1})+d_P(y^{-1}x)) \le \phi(0.292)$,
where $\phi$ is defined by \ref{phidef}.  After possibly interchanging
the roles of $y$ and $y^{-1}$, we may therefore assume that
$$d_P(xy)+d_P(yx)
\le \phi(0.292).$$
It follows that
\Equation\label{peau lisse}
\max (d_P(xyx^{-1}y^{-1}),d_P(x^{-1}y^{-1}xy) )\le \phi(0.292).
\EndEquation

We claim that at least one of the subgroups $\langle
xyx^{-1}y^{-1},yx^{-1}y^{-1}x\rangle$, $\langle x^2,yx^2y^{-1}\rangle$
or $\langle y^2, x^{-1}y^2x\rangle$ has infinite index in $\Gamma$.
If $H_1(M;\QQ)\ne0$, it is immediate that
$\langle xyx^{-1}y^{-1},yx^{-1}y^{-1}x\rangle$ has infinite index.  If $M$ is
closed and contains a semi-fiber $F$, then the image of the inclusion
homomorphism $\pi_1(F)\to\pi_1(M)$ is a normal subgroup $N$ of
$\pi_1(N)$, and $D=\pi_1(M)/N$ is an infinite dihedral group. Hence
the commutator subgroup $D'$ of $D$ is infinite cyclic. If the images
$\bar x$ and $\bar y$ of $x$ and $y$ in $D$ belong to $D'$, then $\bar
x$ and $\bar y$ commute; thus in this case $\langle
xyx^{-1}y^{-1},yx^{-1}y^{-1}x\rangle$ is contained in $N$, and
therefore has infinite index in $\Gamma$. If $\bar x$ does not belong
to $D'$ then $\bar x$ has order $2$ in $D$ and hence $x^2\in N$; thus
$\langle x^2,yx^2y^{-1}\rangle$ is contained in $N$, and therefore has
infinite index in $\Gamma$.  Similarly, if $\bar y\notin D'$, then
$\langle y^2, x^{-1}y^2x\rangle$ has infinite index in $\Gamma$.

Since $M$ is a hyperbolic $3$-manifold, the manifold-with-boundary
obtained from $M$ by removing a standard open cusp neighborhood for
each $\ZZ\times\ZZ$-cusp satisfies the hypothesis of
\cite{JS}*{Theorem VI.4.1}. It therefore follows from the latter
theorem that every two-generator subgroup of infinite index in
$\Gamma=\pi_1(M)$ is either free of rank at most $2$ or free abelian
of rank $2$. Hence at least one of the subgroups $\langle
xyx^{-1}y^{-1},yx^{-1}y^{-1}x\rangle$, $\langle x^2,yx^2y^{-1}\rangle$
or $\langle y^2, x^{-1}y^2x\rangle$ is either free of rank at most $2$
or free abelian of rank $2$.

We shall now assume that $x$ and $y$ do not commute, and deduce a
contradiction.  It follows from Assertion (\ref{langevin}) of Corollary
\ref{commutators} that the elements $ xyx^{-1}y^{-1}$ and
$yx^{-1}y^{-1}x$ do not commute. On the other hand, it follows from 
Assertion (\ref{lauvergnat}) of Corollary \ref{commutators}
that $x^2$ does not commute with $y$, and it
therefore follows from the second assertion of Lemma \ref{commute}  that
$x^2$ does not commute with $yx^2y^{-1}$ or with $y^{-1}x^2y$. Thus the
subgroups $\langle xyx^{-1}y^{-1},yx^{-1}y^{-1}x\rangle$, $\langle
x^2,yx^2y^{-1}\rangle$ or $\langle y^2, x^{-1}y^2x\rangle$ are all
non-abelian. Hence at least one of these subgroups
is free of rank $2$;
that is, at least one of the pairs $( xyx^{-1}y^{-1},yx^{-1}y^{-1}x)$,
$( x^2,yx^2y^{-1})$ or $( y^2, x^{-1}y^2x)$ is independent.

If $xyx^{-1}y^{-1}$ and $yx^{-1}y^{-1}x$ are independent, it follows from 
\cite{acs}*{Theorem 4.1} that
\Equation\label{just kidding}\frac1{1+\exp d_P(xyx^{-1}y^{-1})}+\frac1{1+\exp d_P(yx^{-1}y^{-1}x)}
\le\frac12.
\EndEquation
On the other hand, by (\ref{peau lisse}) we have 
$$
\frac1{1+\exp d_P(xyx^{-1}y^{-1})}+\frac1{1+\exp d_P(yx^{-1}y^{-1}x)}
\ge\frac2{1+\exp\phi(0.292)}=0.5009\ldots,
$$
which contradicts (\ref{just kidding}). 

Now suppose that $x^2$ and $yx^2y^{-1}$ are independent. We have
$d_P(x^2)\le2d_P(x)\le2\cdot0.292$ and
$d_P(yx^2y^{-1})\le2d_P(x)+2d_P(y)\le4\cdot0.292$.  From
\cite{acs}*{Theorem 4.1} we find that
$$\begin{aligned}
\frac12&\ge\frac1{1+\exp d_P(x^2)}+\frac1{1+\exp d_P(yx^2y^{-1})}\\
&\ge\frac1{1+\exp(2\cdot0.292)}+\frac1{1+\exp(4\cdot0.292)}\\
&=0.595\ldots,\end{aligned}
$$
a contradiction. We obtain a contradiction in the same way if $y^2$
and $x^{-1}y^2x$ are independent. This completes the proof of the
proposition in the special case where (\ref{squares still have more
  fun}) holds.

We now turn to the general case, in which the inequality \ref{squares
still have more fun} is not assumed to hold. Since $x$ and $y$ are
loxodromic, the quantities $d_P(x^n)$ and $d_P(y^n)$ tend to
$\infty$ with $n$. Hence there is a largest integer $n_1$ such that
$d_P(x^{n_1})<0.292$, and there is a largest integer $n_2$ such
that $d_P(y^{n_2})<0.292$. If we set $x_0=x^{n_1}$ and
$y_0=y^{n_2}$, it follows that $d_P(x_0^2)\ge0.292$ and that
$d_P(y_0^2)\ge0.292$. We may now apply the special case of the
proposition that has already been proved, with $x_0$ and $y_0$ in the
roles of $x$ and $y$, to deduce that $x_0$ and $y_0$
commute.  It then follows that $x$ and $y$ commute.
\EndProof

\Proof[Proof of Theorem \ref{main theorem or wall theorem}] The second
assertion is Proposition \ref{off the wall}. In proving the first
assertion we may assume that $H_1(M;\QQ)=0$ and that $M$ is not a
closed manifold containing a semi-fiber. The condition $H_1(M;\QQ)=0$
implies that $M$ is closed and not fibered. The proof of the first
assertion is thus reduced to the case where $M$ is closed and contains
an incompressible surface $F$ which is not a fiber or a semi-fiber. In
this case the result follows immediately from assertions (\ref{won})
and (\ref{tree}) of Proposition \ref{sanity clause}.
\EndProof


\end{document}